\documentclass[12pt]{article}
\usepackage{epsfig}
\usepackage{amsmath}
\usepackage{amssymb}

\newcommand{\nit}{\noindent}
\newcommand{\no}{\nonumber}
\newcommand{\be}{\begin{equation}}
\newcommand{\ee}{\end{equation}}
\newcommand{\br}{\begin{eqnarray}}
\newcommand{\er}{\end{eqnarray}}

\newcommand{\lam}{\mbox{$\lambda$}}

\newcommand{\abs}[1]{\lvert #1 \rvert}
\newcommand{\norm}[1]{\lVert #1 \rVert}

\newcommand{\avg}[1]{\langle #1 \rangle}
\newtheorem{theo}{Theorem}[section]

\newtheorem{prop}{Proposition}[section]

\newtheorem{rmk}{Remark}[section]

\begin{document}
\title{Min-Max Variational Principle and
Front Speeds in Random Shear Flows}

\author{James Nolen\thanks{Department of Mathematics, University of Texas at Austin,
Austin, TX 78712 (jnolen@math.utexas.edu).}
\and Jack Xin\thanks{Department of Mathematics and ICES (Institute of Computational
Engineering and Sciences), University of Texas at Austin, Austin, TX 78712
(jxin@math.utexas.edu).}}

\date{}
\setcounter{page}{1}
\setcounter{section}{0}
\maketitle
\begin{abstract}
Speed ensemble of bistable (combustion) fronts in mean zero
stationary Gaussian shear flows
inside two and three dimensional channels
is studied with a min-max variational principle.
In the small root mean square regime of shear flows,
a new class of multi-scale test functions are found to yield speed
asymptotics. The quadratic speed enhancement law holds with probability
arbitrarily close to one under the almost sure continuity (dimension two) and mean
square H\"older regularity (dimension three) of the shear flows.
Remarks are made on the conditions for
the linear growth of front speed expectation
in the large root mean square regime.

\end{abstract}

\newpage

\section{Introduction}
\setcounter{equation}{0}
We consider propagation speeds of reaction-diffusion fronts in
random shear flows. The model equation is:
\be
u_t =\Delta_{x} u +\delta B \cdot \nabla_{x} u+ f(u), \label{e0}
\ee
where $x=(x_1,\tilde x) \in D=R^1 \times \Omega$, $\Omega = [0,L]^{n-1}$,
$n=2,3$; $\delta > 0$ a scaling parameter;
the vector field $B=(b(\tilde x,\omega),0)$,
$b(\tilde x,\omega)$ is a scalar stationary Gaussian field with
zero ensemble mean and almost surely continuous sample paths;
$f$ is either one of the following:
\medskip

\nit (1) the combustion nonlinearity with ignition temperature cutoff:
$f(u) = 0$, $ u \in [0,\mu]$, $\mu \in (0,1)$;
$f(u) > 0$, $u \in (\mu,1)$, $f(1)=0$, $f\in C^2([0,1]$;
\medskip

\nit (2) the bistable nonlinearity:
$f(u)=u(1-u)(u-\mu )$, $\mu \in (0,1/2)$.
\medskip

The boundary condition is zero Neumann at $R^1 \times \partial \Omega$.
The initial data $u_0$ belongs to the set $I_s$. The set
$I_s$ for bistable $f$
consists of bounded continuous functions with
limits one and zero at $x_{1} \sim \pm \infty$ respectively.
For combustion $f$, one requires also that $u_0$ decay to zero
exponentially at $-\infty$.  For each realization of $B$,
there exists a traveling front solution of the form
$u=\Phi(x_1 +c t, \tilde x)$, with unique
speed $c$ and profile $\Phi$ \cite{BN1}. For $u_0 \in I_s$, $u(x,t)$ converges to a
traveling front at large times \cite{Roque}.

We are interested in the asymptotic behavior of $c=c(\delta)$
when $\delta$ is small or large. Such speed asymptotics have been recently
studied for both deterministic fronts \cite{B1,Const,HPS,Kh,KR,NX1,PX} and
KPP fronts in random shears \cite{NX2,Xin3}. See also \cite{Xin2} for
related applied science literature.
The bistable (combustion) front speed $c$ has a min-max variational
characterization \cite{Ham,HPS}.
Following the notation of \cite{HPS}, define the functional:
\be
\psi(v) = \psi(v(x)) \equiv {L\, v + \, f(v) \over \partial_{x_1} v} \equiv
\frac{\Delta v + \delta b(\tilde x)
\partial_{x_1} v + f(v)}{\partial_{x_1} v}. \label{v1}
\ee
The strong form of the
min-max variational formula of \cite{HPS} for $c(\delta)$ is:
\be
\sup_{v \in K} \inf_{x \in D} \psi(v(x)) = c(\delta) =
\inf_{v \in K} \sup_{x \in D} \psi(v(x)) \label{v2}
\ee
where $K$ is the set of admissible functions:
\[ K = \{v \in C^{2}(D)| \; \partial_{x_1}\, v > 0,\,
0 < v(x) < 1, v \in I_s\}. \]

The variational principle (\ref{v2}) is a powerful tool for getting
tight bounds on $c$, if one can construct an admissible
test function $v$ that approximates
well the exact traveling front solution $\Phi$. This is the case when $\delta $
is small and the test function is a perturbation of the known traveling
front solution when $\delta = 0$. In \cite{HPS}, a form of such test function
is proposed for deterministic shear flow, and the quadratic
speed enhancement law (\cite{PX} and references therein)
is recovered if $\delta $ is small enough. To handle an ensemble of
shear flows in the random case, one must know which norm of $b$ controls
the smallness of $\delta $ by analyzing an admissible test function $v$
and the functional $\psi(v)$.
A new form of multi-scale test function is found for this purpose, and
it cures a delicate divergence problem when perturbing the traveling front
at $\delta =0$. It turns out that the smallness of $\delta$ for the quadratic
speed enhancement depends on the maximum (H\"older) norm of $b$ if $n=2 \,(3)$.

Our main results are that the quadratic speed enhancement law holds
with probability arbitrarily close to one, under the almost sure
continuity ($n=2$) or mean square H\"older regularity ($n=3$)
conditions of the shear.
As $b$ can be arbitrarily large (though with small probability),
our admissible test functions are not valid for the entire ensemble.
We are unable to discuss the expectation of front speed. This is in contrast to
the Kolmogorov-Petrovsky-Piskunov (KPP) case (e.g. when $f(u)=u(1-u)$),
where there is less restriction on test functions and the minimal speed
variational principle is minimization over functions on the cross section
$\Omega$ \cite{NX2}.

The paper is organized as follows. In section 2, the new
multi-scale test function is presented and used to estimate the front speed and
identify the smallness of $\delta$ in terms of proper norms of $b$.
In section 3, probabilistic estimates based on Borell's inequality and
Karhunen-Lo\`eve expansion of mean zero Gaussian fields help to control
the maximum (H\"older) norms of $b$, and lead to the quadratic speed
enhancement laws. Remarks are made on the linear front speed growth law
when $\delta $ is large.

\section{Deterministic Front Speed Asymptotics}
\setcounter{equation}{0}
Let us consider $b$ to be a given H\"older continuous deterministic
function, and assume that
$<b> \equiv  {1\over |\Omega |}\int_{\Omega} \, b(y)\, dy = 0$.
Otherwise, the integral mean contributes
a linear term $\delta <b>$ to the front speed. Let $U=U(x_1 + c_0 t)$ be
the traveling fronts when $\delta =0$, satisfying the equation:
\[ U'' - c_0 U' + f(U) =0,\;\; \]
$U(-\infty)=0$, $U(+\infty)=1$, $U' > 0$, $U$ approaches zero and one
at exponential rates.

In the small $\delta $ regime, define the test function
\be
v(x) = U(\xi) + \delta^2 \tilde u(\xi,\tilde x), \label{v3}
\ee
where the variable:
\be \xi = (1 + \alpha \delta^2)x_1 + \delta \chi, \label{v4}
\ee
with $\alpha $ a constant to be determined and
$\chi = \chi(\tilde x)$ solution of:
\[
-\Delta_{\tilde x} \chi = b, \;\; \tilde x \in \Omega,
\]
subject to zero Neumann boundary condition at $\partial \Omega$.
We normalize $\chi$ by setting $\int_\Omega \chi\,dy = 0$.
The test function (\ref{v3}) reduces to the one in \cite{HPS} if $\alpha =0$.
However, we shall see that $\alpha \not = 0$ is essential in
suppressing certain divergence arising in the evaluation of $\psi (v)$,
similar to curing secular growth in multiple-scale perturbation expansions.

For H\"older continuous $b(\tilde x)$, $\chi$ is $C^{2+p}$, $p \in (0,1)$.
Assume that $\tilde u$ is $C^{2}$, and decays sufficiently fast
at $\xi $ infinities. A straightforward computation shows that:
\br
L v + f(v) &=& U'' + f(U) \no \\
& & + \delta^2\left( 2\, \alpha\, U'' + \Delta_{\xi,\tilde x} \tilde u
+ U''\abs{\nabla \chi}^2 + f'(U) \tilde u \right) \no \\
& & + \delta^3 (b \, \alpha \, U'  + 2 \nabla_{\tilde x} \tilde u_\xi \cdot \nabla \chi)\no \\
& & + \delta^4 (\alpha^2\, U'' + 2 \, \alpha\, \tilde u_{\xi \xi} + \Delta_\xi \,
\tilde u \, \abs{\nabla \chi}^2) \no \\
& & + f(v) - f(U) - f'(U)\,\delta^2 \, \tilde u  \no \\
& & + \delta^5\, \alpha \, b \, \tilde u_\xi  + \delta^6\, \alpha^2 \,
\tilde u_{\xi \xi},
\label{v5}
\er
where $f(v) - f(U) - f'(U)\delta^2 \tilde u = O(\delta^4 \tilde u^2)$ and
the $O(\delta)$ terms cancel by choice of $-\Delta \chi = b$.
On the other hand, up to an undetermined constant $\gamma$, we have:
\br
(c_0 + \delta^2 \gamma)\; \partial_{x_1}\, v \,
&=& c_0 U' + \delta^2(c_0\alpha U' + c_0 \tilde u_\xi + \gamma U') \no \\
& & + \delta^4(\gamma \alpha U' + c_0 \alpha \tilde u_\xi + \gamma \tilde u_\xi) \no \\
& & + \delta^6(\alpha \gamma \tilde u_\xi). \; \label{v6}
\er
All terms in (\ref{v5}) and (\ref{v6}) are evaluated
at $(\xi,\tilde x)$. Let us choose $\tilde u$ to solve the equation:
\be
\bar L \, \tilde u \, = \, \Delta_{\xi,\tilde x} \, \tilde u -
c_0 \, \tilde  u_\xi + f'(U) \, \tilde u = -\abs{\nabla \chi}^2\, U'' -
2\, \alpha\, U'' + \gamma\, U' + c_0\, \alpha \, U', \label{e2}
\ee
subject to zero Neumann boundary condition at $\partial \Omega $,
and exponential decay at $\xi $ infinities.
The solvability of (\ref{e2}) will be discussed later.
Then we would have
\be
L v + f(v) = (\partial_{x_1} v) (c_0 + \delta^2 \gamma) + R_1. \label{v8}
\ee
The remainder $R_1$ is:
\be
R_1 = \delta^3 A + \delta^4 B + \delta^5 D + \delta^6 E + O(\delta^4
\tilde u^2), \label{e34}
\ee
where
\br
A &=& b \, \alpha U'  + 2 \nabla_{\tilde x}\, \tilde u_\xi \cdot \,
\nabla \chi \no \\
B &=& \alpha^2\, U'' + 2 \alpha \,\tilde u_{\xi \xi} + \,\Delta_\xi \,
\tilde u \, \abs{\nabla \chi}^2 - \gamma \, \alpha \, U' - c_0 \, \alpha \,
\tilde u_\xi - \gamma \, \tilde u_\xi \no \\
D &=& \alpha b \, \tilde u_\xi \no \\
E &=& \alpha^2 \, \tilde u_{\xi \xi} - \alpha \, \gamma \, \tilde u_\xi. \label{v9}
\er

The constant $\gamma$ is determined by a solvability condition for (\ref{e2}).
The right hand side must be orthogonal to the
function $U'(\xi)\, e^{-c_0 \xi}$
which spans the kernel of the adjoint operator
$\bar L^* v = \Delta_{\xi, \tilde x} v + c_0 v_\xi + f'(U) v$.
Thus,
\br
\gamma \int_R (U')^2 e^{-c_0 \xi}\,d\xi &=& \frac{1}{\abs{\Omega}}\int_D \left((\abs{\nabla \chi}^2 + 2\alpha ) U''U' - c_0 \alpha (U')^2\right)  e^{-c_0 \xi} d\xi\, d \tilde x \no \\
& = & \frac{c_0}{2} \avg{\abs{\nabla \chi}^2 + 2 \alpha }
\int_R (U')^2 e^{-c_0 \xi}\,d\xi - c_0 \alpha \int_R (U')^2
e^{-c_0 \xi}\,d\xi \no
\er
where the braket denotes the integral average over $\Omega$. We have
performed integration by parts once, the boundary terms decay to zero
and the integral $\int_R (U')^2 e^{-c_0 \xi}\,d\xi $ converges for both
bistable and combustion type nonlinearities.
It follows that equation (\ref{e2}) is solvable when
\be
\gamma = \frac{c_0}{2\abs{\Omega}} \int_\Omega \abs{\nabla \chi}^2 \,d\tilde x
 = {c_0 \over 2} \avg{\abs{\nabla \chi}^2}, \label{v10}
\ee
regardless of the choice of $\alpha$.
For such $\gamma$,
we may choose $\alpha = -\frac{1}{2} \avg{\abs{\nabla \chi}^2}$
so $\gamma + c_0 \, \alpha = 0$ and the $\tilde u$ equation becomes:
\be
\bar L \tilde u = \Delta_{\xi,\tilde x} \tilde u -
c_0 \tilde u_\xi + f'(U) \tilde u =
\left(\avg{\abs{\nabla \chi}^2} -\abs{\nabla \chi}^2\right) U''. \label{v11}
\ee
Note that the right hand side of (\ref{v11}) has zero integral average
over $\Omega$.

The solution $\tilde u$ of (\ref{v11})
can be expanded in
eigenfunctions of the Laplacian $\Delta_{\tilde x}$:
\br
\tilde u(\xi, \tilde x) & = & \sum_{j \in Z^{n-1}_{+}\cup \{0\}}\,
u_j(\xi) \phi_j(\tilde x), \no \\
\abs{\nabla_{\tilde x} \chi}^2 & = & \sum_{j \in Z^{n-1}_{+}\cup \{0\}}
\, a_j \phi_j(\tilde x), \no
\er
where $Z^{n-1}_{+}$ denotes nonnegative integer vectors in $R^{n-1}$
with at least one positive component;
$-\Delta_{\tilde x} \phi_j = \lambda_j \phi_j$ with
zero Neumann boundary condition on $\partial \, \Omega$. The eigenfunctions
$\phi_j$ are normalized so that $\|\phi_j\|_{L^{\infty}} =1$, and
the eigenvalues $\lambda_j = O(|j|^{2})$ for large $|j|$.
Clearly, $\phi_0 = 1$ and $\lambda_0 = 0$, $a_0 = \avg{\abs{\nabla \chi}^2}$.
The functions $u_j$ solve the equations
\be
L_j u_j = u_j'' - \lambda_j u_j - c_0 u_j' + f'(U) u_j =
-(a_j -\delta_{j,0} \, a_0)\, U''.
\label{e1}
\ee

The $j=0$ equation has zero right hand side, and we take $u_0 \equiv 0$.
If $\alpha$ were zero, equation (\ref{e2}) implies that
the right hand side of the $j=0$ equation would
be $- a_0\, U'' + \gamma \, U'$. The general solution (zeroth-mode) for
$u_0$ is $-{a_{0}\over 2} \xi \, U'(\xi) + {\rm const.}\, U'(\xi)$.
The linear factor $\xi$ would
render $v(x)$ unable to stay in the set $K$ due to
divergence at large $\xi$ (and $x_1$).

Now consider the equations for $j \not = 0$.
As shown in \cite{Xin1}, the operators $L_j$
are invertible in part due to $\lambda_j > 0$ for $j \not = 0$.
As in Lemma 2.3 of \cite{Xin1}, we have the regularity estimate:
\be
\norm{u_j}_{C^2(R)} \leq C_1 \abs{a_j}, \;\;\;\; j\not = 0, \label{e16}
\ee
where the constant $C_1$ is independent of $j$ and
depends only on $U$ and its derivatives.

Next we show that for $j\not = 0$, the ratio $\frac{u_j}{U'}$ is bounded by
$C_2 \, |a_j|$, where $C_2$ is a positive constant independent of $j$.
It suffices to prove a uniform bound if the $a_j$ factor is one.
Consider the function $w_j = \beta U' - u_j$ for some positive
constant $\beta$. For some $r_0 > 0$ sufficiently large,
$f'(U) \leq 0$ whenever $\abs{\xi} \geq r_0$. Choose $\beta$ as
\be
\beta \geq \beta_1 \equiv \frac{C_1}{\min_{\abs{\xi}\leq r_0} U'(\xi)},
\label{v12}
\ee
so that by (\ref{e16}), $w_j > 0$ in the region $\abs{\xi} \leq r_0$.
By equation (\ref{e1}) for $u_j$ and the fact that
$U'$ solves $L_j U' = - \lambda_j U' < 0$, we have
\be
L_j w_j = - \beta \lambda_j U' + U''. \label{e22}
\ee
Using the fact that $U''$ is bounded by a constant multiple of $U'$, and
that $\inf_{j \not = 0}\, \lambda_j \, > \, 0$,
we may further increase $\beta$ if necessary to ensure
that $L_j w_j < 0$ in the region $\abs \xi > r_0$.
Maximum principle implies that $w > 0 $ for all $\xi \in R$.
Repeating the argument for $-u_j$ and taking into account the $a_j$ factor give:
\be
\frac{\abs{u_j}}{U'} \leq C_3 \abs{a_j}, \label{e35}
\ee
where the constant $C_3$ depends only on $U$,
and (\ref{e35}) holds uniformly
in $j \not = 0$. Inequality (\ref{e35}) says that there is no resonance
when inverting the operator $L_j$ to find $u_j$, for $j \not = 0$, in other words
$u_{j}$ decays same as $U'$ at infinities.

To improve the estimate above, we move the term $f'(U)\,u_j$ to
the right hand side of (\ref{e1}) which
is then bounded by $C'_{3} |a_{j}| U'$, for $C'_{3}$ independent of $j$.
The left hand side operator becomes $u_j'' - \lambda_j \, u_j - c_0\,  u_j'$.
An explicit formula can be written for $u_j$, and the asymptotics of
$\lambda_j$ imply that
\be
{\abs{u^{(i)}_j(\xi,\tilde x)} \over U'(\xi)} \leq C_4
{\abs{a_j}\over (1+|j|^2)}, \label{e18}
\ee
for a constant $C_4$ depending only on $U$, where $i=0,1,2$ denotes
the order of $\xi$ derivatives.

Suppose that $b$ is
H\"older continuous, then Schauder estimates give \cite{GT}:
\be
\norm{\abs{\nabla \chi}^2 }_{1+p} \leq C \, \norm{b}^2_{p},
\label{v14a}
\ee
for some $p\in (0,1)$, $C=C(\Omega)$, $\| \cdot \|_{p}$ is the standard
H\"older norm.
For rectangular cross section $\Omega $ (dimension $n-1$),
the $\phi_j$'s are trigonometric functions. Then
\be
 \abs{a_j} \leq C_5 \;
\norm{b}^2_{p} \; (1+|j|)^{-(1+p)},\;\; p \in (0,1), \label{e7}
\ee
for some constant $C_5$ \cite{Zyg}. Combining (\ref{e18}), (\ref{v14a}), and
(\ref{e7}), we see that the eigenfunction expansion of $\tilde u$
\be
\tilde u(\xi,\tilde x) = \sum_{j\in Z^{n-1}_{+}} \, u_j(\xi)\,
\phi_j(\tilde x) \label{e19}
\ee
converges uniformly in $(\xi,\tilde x)$. Moreover,
\br
{|\tilde u^{(i)}_{\xi}(\xi,\tilde x)| \over U'(\xi)} & \leq &
\sum_{j\in Z^{n-1}_{+}}\,
|u^{(i)}_{j}|(\xi)\,|\phi_j|(\tilde x)/U' \no \\
& \leq & C_4\, C_5\,\norm{b}^2_{p} \; \sum_{j\in Z^{n-1}_{+}}
\, (1+|j|)^{-(3+p)} =C_6\,\norm{b}^2_{p}, \label{v13}
\er
if $n=2,3$. The mixed derivative term is bounded as:
\br
{|\nabla_{\tilde x}\,\tilde u_{\xi}(\xi,\tilde x)|\over U'(\xi)} & \leq &
\sum_{j\in Z^{n-1}_{+}}\, {|u^{(1)}_{j}|(\xi)\over U'}\, |j| \,
 |\phi_j|(\tilde x) \no \\
& \leq & C_5\, \norm{b}^{2}_{p} \;
\sum_{j\in Z^{n-1}_{+}} \, (1+|j|)^{-(2+p)}
 = C_7\,\norm{b}^{2}_{p}, \label{v14}
\er
if $n=2,3$. In case $n=2$, $\Omega$ is an interval, then
$b \in L^{\infty}(\Omega)$ suffices, because
$\nabla_{\tilde x}\, \chi $ is Lipschitz, so
$|a_j| \leq O(\|b\|_{\infty}^{2} |j|^{-p})$
for some $p \in (0,1)$ which replaces estimate (\ref{e7}).
The exponent in (\ref{v14})
goes down to $1+p$, yet enough for convergence on $j \in Z^{1}_{+}$.
As a result, for $n=2$, $i=0,1,2$, the estimates:
\be
{|\tilde u^{(i)}_{\xi}(\xi,\tilde x)| \over U'(\xi)} \leq
C_8\, \|b\|^{2}_{\infty}, \;\;\;
{|\nabla_{\tilde x}\,\tilde u_{\xi}(\xi,\tilde x)|\over U'(\xi)}  \leq
C_9\, \|b\|^{2}_{\infty}, \label{v15}
\ee
hold.

As a consequence of these estimates, let us verify the admissibility of
test function $v$. Clearly, $v$ approaches zero (one) exponentially at $x_1 =
-\infty$ $ (+\infty)$. As
\[v_{x_1}=(1 + \delta^2 \alpha) U' +
\delta^2 (1 + \delta^2 \alpha) \tilde u_\xi, \]
we have for $n=3$:
\be
v_{x_1} \geq \frac{1}{2} U' - \frac{1}{2} \delta^2 \abs{\tilde u_\xi}
 \geq  \left [\frac{1}{2} - \delta^2 \, \frac{C_6}{2} \, \norm{b}^2_{p}\,\right ]\, U'
 \geq {1\over 4}\, U' > 0, \label{v16}
\ee
if:
\be
 \delta^2 \leq \min \left (\frac{1}{2 \abs{\alpha}},
{1\over 2}\, (C_{6}\,\norm{b}^2_{p})^{-1}\right ). \label{v17}
\ee
For $n=2$, in view of (\ref{v15}), $v_{x_1}\, \geq \, {1\over 4} \, U'$ if
\be
\delta^2 \leq \min \left (\frac{1}{2 \abs{\alpha}},
{1\over 2}\, (C_{8}\,\norm{b}^2_{\infty})^{-1}\right ). \label{v18}
\ee
It follows from the monotonicity of $v$ in $x_1$ and its
asymptotics near $x_1$ infinities that $0 < v < 1$ for all $x_1,\tilde x$.

With $v$ being an admissible test function, we see from (\ref{v8}) that
\be
\psi(v(x)) = \frac{L v + f(v)}{\partial_{x_1} v} =
c_0 + \delta^2 \gamma + \frac{R_1}{\partial_{x_1} v}, \label{v19}
\ee
with $R_1$ defined by (\ref{e34}). By (\ref{v17}) or (\ref{v18}),
 $\partial_{x_1} v \geq \frac{1}{4} U'$. Also,
each term in $R_1$ is bounded by a multiple of $U'$.
It follows that for $\delta$ satisfying (\ref{v17}) or (\ref{v18}),
$
\abs{\frac{R_1}{\partial_{x_1} v}} \leq \abs{4\frac{R_1}{U'}} \leq \, C_{10}\, \delta^3,
$
where $C_{10} = C_{10} (U, \| b\|_{p}), p \in (0,1)$, if $n=3$,
$C_{10} = C_{10} (U, \|b\|_{\infty})$ if $n=2$.
We have shown:
\begin{prop}
Let $n=2$ or $3$, and let the shear flow profile
$b(\tilde x)$ be H\"older continuous with
exponent $p \in (0,1)$ over rectangular domain $\Omega \subset R^{n-1}$.
There is a positive constant $\delta_0$ depending only on $U$ and
the H\"older (maximum) norm of $b$ for $n=3$ ($n=2$) such that
if $\delta \leq \delta_0$:
\be
c(\delta) = c_0 + \frac{c_0 \delta^2}{2 \abs{\Omega}}
\int_\Omega \abs{\nabla \chi}^2 \,dy + O(\delta^3).
\ee
A specific form of $\delta_0$ follows from (\ref{v17}) and (\ref{v18}).
\end{prop}

At large $\delta$, numerical
evidence suggested that $\lim_{\delta\to \infty}\, c(\delta)/\delta$ exists for
bistable and combustion nonlinearities \cite{NX1}.
Such a limit holds for the KPP nonlinearity \cite{B1}.

\section{Random Front Speed Asymptotics}
\setcounter{equation}{0}
Let us consider $b=b(\tilde x,\omega)$ as a stationary mean zero Gaussian
field with almost surely continuous sample paths in dimension one and
two for $n=2,\, 3$ respectively. We are interested
in the restriction of $b$ over $\Omega$. In order to apply results of the
previous section on each realization of $b$, let us write $b =\bar{b} + b_{1}$,
where $\bar{b} = |\Omega |^{-1}\, \int_{\Omega}\, b(\tilde x,\omega)\, d\tilde x$.
Correspondingly, $c=c(\delta,\omega) = c_0 - \delta \bar{b} + c_2 (\omega)$,
where $c_{2}$ depends on $b_{1}$. As the size of $\delta_0$ depends on either the
maximum norm or the H\"older norm of $b_1$ which is an unbounded random variable,
there is no uniform way to choose a $\delta_0$ for all realizations.
On the other hand, the probability of the occurance of very
large maximum (H\"older) norm of $b_1$ is often small, so the random speed
asymptotics may hold with probibility arbitrarily close to one.
Let us consider $n=2$ first.

\begin{theo}[Two dimensional channel]
Let $D = R^1 \times [0,L]$, $\Omega = [0,L]$, and $b(\tilde x)$ be the restriction on
$\Omega$ of a mean zero stationary, Gaussian random process with almost surely
continuous sample paths. Let $f$ be a bistable or combustion nonlinearity.
Then for each small $\epsilon \in (0,1/4)$ and $q \in (0,1)$, there is a deterministic constant $\delta_0 = \delta_0 (\epsilon, q)$
such that if $\delta \in (0,\delta_0)$,
\[ Prob \left \{ \left \lvert c(\delta,\omega) - c_0 + \delta \bar{b} - \delta^{2}\,
{1\over 2}\avg{\left (\int_{0}^{\tilde x}\, b_{1}(y)\, dy \right )^{2}} \right \rvert
\geq \kappa \delta^{3-q}  \right \} < \epsilon. \]
where the constant $\kappa>0$ is independent of $\delta$, $\epsilon$, and $q$. As $\epsilon \to 0$, the constant $\delta_0$ can be chosen to satisfy
\be
\delta_0(\epsilon,q) \geq C  \abs{\log(\epsilon)}^{-\frac{2}{q}}
\ee
\end{theo}

\nit {\bf Proof:} Recall the inequality (Lemma 3.1,
p. 62, \cite{Adler2}):
\be
E[\sup_{\Omega}\, b(\tilde x)] \leq E\sup_{\Omega}\, |b(\tilde x)|
\leq E |b(0)| + 2 E [\sup_{\Omega}\, b(\tilde x)],
\ee
and the Borell inequality on mean zero Gaussian process with
almost surely continuous sample paths (Theorem 2.1, p. 43, \cite{Adler2}):
\br
& & \mu \equiv E[ \sup_{\Omega}\, b(\tilde x) ] < \infty \no \\
& & P(|\sup_{\Omega}\, b(\tilde x) -\mu | > \lambda ) \leq 2
e^{-\lambda^2/(2\sigma^2)}, \label{b1}
\er
where $\sigma^2 = E[b^2]$. With the choice of $\lambda = \sqrt{ - 2 \sigma^2 \log(\epsilon / 4)}$, inequality (\ref{b1}) implies that
\br
P( \norm{b}_\infty > \lambda + \mu) &\leq&  P( \sup_\Omega b(\tilde x) - \mu > \lambda) + P( - \inf_\Omega b(\tilde x) - \mu > \lambda) \no \\
& \leq & 4 e^{-\lambda^2/(2\sigma^2)} \leq \epsilon .
\er
We now see that for $\epsilon>0$ sufficiently small, $\lambda > 1$ and there is a constant $C_{11}$ independent of $\epsilon$ such that if
\be
\delta < \delta_0(\epsilon) = C_{11} \left( \frac{1}{\mu + \lambda}\right)^\frac{4}{q} \leq C_{11} \frac{1}{\mu + \lambda} \label{d1}
\ee
then (\ref{v18}) holds with probability at least $1-\epsilon$. Therefore, by Proposition 2.1, the front speed asymptotics hold with probability $1-\epsilon$.  Using (\ref{e34}), (\ref{v9}), and (\ref{v15}) we see that the remainder in the expansion satisfies
\br
\left \lvert \frac{R_1}{\partial_{x_1} v} \right \rvert \leq C_{12} \delta^3 \norm{b}^4_\infty & \leq & C_{12} \delta^3 (\mu + \lambda)^4 \no \\
& \leq & C_{13} \delta^{3-q} \delta^q \left( \mu + \lambda \right)^4  \no \\
& \leq & C_{13} \delta^{3-q} \delta_0^q \left( \mu + \lambda \right)^4  \no \\
& \leq & C_{14} \delta^{3-q}  \no
\er
with probability $1-\epsilon$. The estimate on $\delta_0(\epsilon)$ follows from (\ref{d1}).

\begin{rmk}
An example of such $b$ process is the mean zero
stationary Ornstein-Uhlenbeck (O-U) process
which is both Gaussian and Markov. The O-U sample paths are almost surely
continuous.
\end{rmk}

The $n=3$ case requires a probabilistic estimate of H\"older norm of
$b$. For a mean zero Gaussian field, this is related to
the structure of the covariance function
$R(t,s)=E[b(t)b(s)]$, $t$, $s \in \Omega$. If
the covariance function is continuous, positive and non-negative definite,
there exists a Gaussian process with this covariance \cite{Adler2}.
The symmetric integral operator:
$\phi \to \int_{\Omega}\, R(t,s)\, \phi(t) \, dt$
generates a complete set of orthonormal eigenfunctions $\phi_j$
on $L^2(\Omega)$ with nonnegative eigenvalues $\lambda_j$, $j=1,2,\cdots$.
Define:
\[ p(u) = \max_{\|s - t\| \leq |u|\sqrt{2}}\, [E |b(s)-b(t)|^2]^{1/2}. \]
Consider the partial sum ($m$ a positive integer):
\be
X^{(m)}(t,\omega)=\sum_{j=1}^{m}\, \sqrt{\lam_j}\,\phi_{j}(t)\theta_{j}(\omega),
\label{b2}
\ee
where $\theta_j$'s are independent unit Gaussian random variables.
The convergence of the partial sum to $b$ is
given by Garsia's Theorem (Theorem 3.3.2, p. 52, \cite{Adler1}).
It says that if $\int_{0}^{1}\, (-\log u)^{1/2}\, dp(u) < \infty$,
then with probability one, $X^{(m)}(t)$ are almost surely equicontinuous and
converge uniformly on $\Omega$. The resulting infinite series is the
celebrated Karhunen-Lo\`eve expansion.
Moreover, the following estimate holds for all $m$:
\br
|X^{(m)}(s) -X^{(m)}(t)| & \leq & 16\sqrt{2}[\log B]^{1/2}\, p(\|s -t\|) \no \\
& & + 32\sqrt{2}\int_{0}^{\|s - t\|} \, (-\log u)^{1/2}\, dp(u), \label{b3}
\er
where $B=B(\omega)$ is a positive random variable, $E[B^2] \leq 32$.

Suppose that $p(u)$ is H\"older continuous with exponent $s \in (0,1)$, then $b$ is almost surely H\"older
continuous with exponent $s$, and the H\"older norm of $b$ is bounded by
 $\alpha_1 [\log B]^{1/2} + \alpha_2$, where $\alpha_1$, $\alpha_2$ are
two positive deterministic constants. Chebyshev's
inequality gives:
\be
{\rm Prob} \, \left( [\log B]^{1/2} \geq \lam \right)=\,{\rm Prob}\, ( B \geq e^{\lam^2} )
\leq E(B^2)/e^{2\lam^2} \leq 32 \, e^{-2\lam^2}. \label{b4}
\ee
This implies that
\br
{\rm Prob}\, \left( \norm{b}_s > \lambda \right) \leq 32 e^{-2 \left(\frac{\lambda - \alpha_2}{\alpha_1} \right)^2} \leq C_{15} e^{-2\lambda^2} \leq \epsilon ,
\er
if $\lambda = \sqrt{\log(\frac{\epsilon}{C_{15}})}$. For $\epsilon >0$ sufficiently small, $\lambda > 1$ and we take
\be
\delta < \delta_0(\epsilon) = C \left( \frac{1}{\lambda} \right)^{\frac{4}{q}} \leq C \frac{1}{\lambda}, \no
\ee
so that (\ref{v17}) holds with probability at least $1-\epsilon$. As before we use (\ref{e34}), (\ref{v9}), and (\ref{v15}) to conclude that the remainder in the expansion satisfies
\br
\left \lvert \frac{R_1}{\partial_{x_1} v} \right \rvert \leq C_{12} \delta^3 \norm{b}^4_s \leq C_{16} \delta^{3-q}  \no
\er
with probability $1-\epsilon$. We now conclude from Proposition 2.1:

\begin{theo}[Three dimensional channel]
Let $D = R^1 \times \Omega$, $\Omega = [0,L]^2$,
and $b$ be a mean zero Gaussian process
such that the function $p(u)$ is H\"older continuous. Then $b$ has almost surely
H\"older continuous sample paths.  For each small $\epsilon \in (0,1/5)$ and $q \in (0,1)$
there is a deterministic constant $\delta_0 = \delta_0 (\epsilon, q)$
such that if $\delta \in (0,\delta_0)$,
\[ Prob \left \{ \left \lvert c(\delta,\omega) - c_0 + \delta \bar{b} - \delta^{2}\,
{1\over 2}\avg{|\nabla_{\tilde x} \chi|^{2}} \right \rvert
\geq \kappa \delta^{3-q}  \right \} < \epsilon, \]
where $\chi$ satisfies: $-\Delta_{\tilde x}\, \chi = b_1$,
subject to zero Neumann boundary condition at $\partial \Omega$. The constant $\kappa>0$ is independent of $\delta$, $\epsilon$, and $q$.
As $\epsilon \to 0$, the constant $\delta_0$ can be chosen to satisfy
\be
\delta_0(\epsilon,q) \geq C \abs{\log(\epsilon)}^{-\frac{2}{q}} \,.
\ee
\end{theo}

\begin{rmk}
Suppose that $\lim_{\delta \to \infty}\, c(\delta,\omega)/\delta$ exists
almost surely and that $E(\|b\|_{\infty})$ is finite. The dominated convergence
theorem implies that

\nit $\lim_{\delta \to \infty}\, E[c]/\delta $ exists and
is finite. This argument is same as for the KPP case, see \cite{NX2} for details.
\end{rmk}

\section{Concluding Remarks}
Bistable and combustion front speeds in mean zero Gaussian random shear
flows have been studied with the min-max variational principle
of \cite{HPS}. The quadratic speed enhancement law is valid with
probability arbitrarily close to one in both two and three dimensional
channels under the almost sure continuity and the
mean square H\"older regularity conditions of the Gaussian shear flows.

It would be interesting to extend results here to the case of general
convex cross section $\Omega$ with smooth boundary, or to the case of
$\Omega$ in dimension higher than two.

\section{Acknowledgements}
J.X would like to thank Professor M. Cranston for helpful communications.
The work is partially supported by NSF grant ITR-0219004.
J.N. is grateful for support through a VIGRE graduate fellowship at UT
Austin.

\bibliographystyle{plain}

\end{document}